\documentclass[a4paper,12pt]{article}
\usepackage{amssymb}
\usepackage{amsthm}
\newtheorem{Lemma}{Lemma}
\newtheorem{Proposition}{Proposition}
\newtheorem{Theorem}{Theorem}
\newtheorem{Remark}{Remark}
\newtheorem{Corollary}{Corollary}

\newtheorem{Notation}{Notation}
\newcommand{\poh}{p_0, \ldots ,p_h}

\newcommand{\OO}{{\mathcal{O}}}

\newcommand{\FF}{{\mathbb{F}}}   
\newcommand{\PP}{\mathbb{P}}
\newcommand{\CC}{\mathbb{C}}

\begin{document}
\title{Osculating spaces to secant varieties}
\author{E.~Ballico, C.~Bocci, E.~Carlini, C.~Fontanari\thanks{This research is 
part of the T.A.S.C.A. project of I.N.d.A.M., supported by P.A.T. (Trento) and 
M.I.U.R. (Italy).} }
\date{}
\maketitle 

\begin{small}
\begin{center}
\textbf{Abstract}
\end{center}
We generalize the classical Terracini's Lemma to higher order osculating 
spaces to secant varieties. As an application, we address with the 
so-called Horace method the case of the $d$-Veronese embedding of the 
projective $3$-space.
\end{small}

\vspace{0.5cm}

\noindent
\textsc{A.M.S. Math. Subject Classification (2000)}: 14N05.


\section{Introduction}

Let $X \subset \PP^r$ be an integral non-degenerate projective variety of 
dimension $n$. The $h$-secant variety $S^h(X)$ of $X$ is the Zariski 
closure of the set of points in $\PP^r$ lying in the span of $h+1$ 
independent points of $X$. The variety $X$ is said to be $h$-defective 
if $S^h(X)$ has dimension strictly less than the expected 
$\min \{ (h+1)(n+1)-1, r \}$. The dimension of $S^h(X)$ can be computed 
via the classical Terracini's Lemma (\cite{Terracini:11}), asserting 
that, for $\poh$ general points on $X$ and $p$ general in their span, 
one has 
$$
T_p(S^h(X)) = < T_{p_0}(X), \ldots, T_{p_h}(X) >.
$$
Let $m$ be a non-negative integer and let $T^m_p(X)$ denote the 
$m$-osculating space to $X$ at $p$ (in particular, notice that 
$T^0_p(X)= \{p\}$ and $T^1_p(X)= T_p(X)$, the usual tangent space). 
Here we are going to prove the following generalization of Terracini's 
Lemma (see \cite{BalFon:04} for a completely different one, computing 
the tangent space to the secant variety of an osculating variety): 

\begin{Theorem}\label{maintheorem}
Fix non-negative integers $h$ and $m$, let $\poh$ be general ponts on $X$ 
and let $p \in S^h(X)$ be general in $< \poh >$. Then 
$$
T^m_p(S^h(X)) = < T^m_{p_0}(X), \ldots, T^m_{p_h}(X) >.
$$
\end{Theorem}

In particular, it turns out that $T^m_p(S^h(X))$ has dimension much less 
than the expected one (in classical terminology, $S^h(X)$ satisfies many 
differential equations of order $m$ at a general point $p$). However, 
this fact is not so surprising: as pointed out in the introductions of 
\cite{Dye1} and \cite{Dye2}, a variety containing a linear space always 
satisfies several differential equations. Theorem~\ref{maintheorem} is 
indeed a consequence of a more general result (see Proposition~\ref{join}) 
to be proved in Section~\ref{terracini}, where the connections with 
Terracini's results on varieties satisfying many differential equations 
(\cite{Terracini:12}) are also discussed. As an application of 
Theorem~\ref{maintheorem}, we can determine the dimension of the 
osculating spaces to the secant varieties of some remarkable projective 
varieties. We point out the following straightforward consequences: 

\begin{Corollary}\emph{(Ciliberto, Cioffi, Miranda, Orecchia \cite{CCMO:03})} 
For any integer $d \ge 1$, let $V_{2,d}$ denote the $d$-Veronese embedding 
of the projective plane $\PP^2$. Fix integers $h \ge 1$ and $1 \le m \le 20$; 
if $h=1$, assume either $d < m$ or $d >  2m-2$; 
if $h=2$, assume either $d < 3m/2$ or $d >  2m-2$;  
if $h=4$, assume either $d < 2m$ or $d >  (5m-2)/2$;   
if $h=5$, assume either $d < 12m/5$ or $d >  (5m-2)/2$;    
if $h=6$, assume either $d < 21m/8$ or $d >  (8m-2)/3$;  
if $h=7$, assume either $d < 48m/17$ or $d >  (17m-2)/6$. 
Then the homogeneous linear system $\OO_{\PP^2}(d)(-mp_0 \ldots -mp_h)$ 
is non-special. In particular, we have 
$$
\dim T^{m-1}_p(S^h(V_{2,d})) = (h+1){{m+1}\choose{2}} - 1.
$$ 
\end{Corollary}

\begin{Corollary}\emph{(Laface \cite{Laface:02})}
For any integer $n \ge 0$, let $\FF_n$ be a Hirzebruch surface, i.e. 
$\FF_n = \PP(\OO_{\PP^1} \otimes \OO_{\PP^1}(n))$, and let $F$, $H$ 
be the two generators of $\mathrm{Pic}(\FF_n)$ such that $F^2=0$, 
$H^2=n$, $F.H=1$. If $m \le 3$ and $(n,a,b,m,h+1) \notin 
\{ (1,0,4,2,5), (1,0,6,3,5), (5,1,4,3,10), (6,0,4,3,11), (n,2e,2,2,2e+n+1), 
(n,e,0,2,r), (n,4e+n+1,2,3,2e+n+1), (n,3e+1,3,3,2e+n+1), (n,3e,3,3,2e+n+1), 
(n,e,1,3,r), (n,e,0,3,r)\}$, then the homogeneous linear system 
$\OO_{\FF_n}(aH+bF)(-mp_0 \ldots -mp_h)$ is non-special. In particular, 
we have
$$
\dim T^{m-1}_p(S^h(\FF_n)) = (h+1){{m+1}\choose{2}} - 1.
$$ 
\end{Corollary}

We are going to prove also the following result:

\begin{Theorem}\label{horacetheorem}
Fix integers $h \ge 1$, $d \ge 4$, and let $V_{3,d}$ 
denote the $d$-Veronese embedding of the projective space $\PP^3$. 
If either $10(h+1)+5d-12 \le {{d+3}\choose{3}}$ or $10(h+1) \ge 
{{d+3}\choose{3}} +5d-12$, then the homogeneous linear system 
$\OO_{\PP^3}(d)(-3p_0 \ldots -3p_h)$ is non-special. In particular, we have 
$$
\dim T^{2}_p(S^h(V_{3,d})) = 10 (h+1) - 1.
$$ 
\end{Theorem}

Theorem~\ref{horacetheorem} is only a special case of a couple of results 
concerning the postulation of double and triple general points in the 
projective space $\PP^3$ (see Proposition~\ref{a12.1} and 
Proposition~\ref{a12.2}) obtained in Section~\ref{horace} via the 
so-called Horace method for zero-dimensional schemes (see Lemma~\ref{a0}).

We work over the complex field $\CC$. 

We are grateful to Stephanie Yang for providing us with the complete 
list of special systems of fat points with multeplicity at most $3$ 
in the projective plane, which is certainly well-known but seems to be 
written nowhere. 

\section{Terracini's Lemma and Laplace equations}~\label{terracini}

Let $X_0, \ldots, X_h \subseteq \PP^r$ be integral non-degenerate 
projective varieties of dimension $n_i$ and let 
$$
J = J(X_0, \ldots, X_h) := \overline{\bigcup_{p_i \in X_i} < \poh >} 
$$
be the join variety of $X_0, \ldots, X_h$. 
Recall that if $p$ is a smooth point of an integral variety 
$X \subset \PP^r$ of dimension $n$ and 
\begin{eqnarray*}
U \subseteq \CC^n &\longrightarrow& \CC^{r+1} \setminus \{ 0 \} \\
t &\longmapsto& p(t)
\end{eqnarray*}
is the lifting of a local parametrization of $X$ centered in $p$, 
then the \emph{$m$-osculating space} $T_p^m(X)$ is the projective subspace 
spanned by the points $[p_I(0)] \in \PP^r$, where $I = (i_1, \ldots, i_n)$ 
is a multi-index such that $\vert I \vert \le m$ and
$p_I = \frac{\partial^{\vert I \vert}p}{\partial t_1^{i_1} \ldots 
\partial t_n^{i_n}}$. With this notation, the following holds:

\begin{Proposition}~\label{join} 
Let $p_i \in X_i$, $i=0, \ldots, h$, be general points and let $p \in J$ 
be general in $< \poh >$. We have
$$
T_p^m(J) = < T_{p_0}^m(X_1), \ldots, T_{p_h}^m(X_h) >.
$$
\end{Proposition}

\proof
We write the proof in the $h=1$ case, the generalization to arbitrary $h$ 
being obvious. Let $U_i\subseteq X_i $ be an analytic neighborhood of $p_i$ 
and denote by $\Phi_i:\mathbb{A}^{n_i}\rightarrow\PP^r$ a local 
parameterization of $X_i$. With this notation, a parameterization of 
$J(X_1,X_2)$ in a neighborhood of $p$ is given by:
\[\begin{array}{cccc}
\Phi: & \mathbb{A}^{n_1}\times\mathbb{A}^{n_2}\times\mathbb{A}^{2} 
& \longrightarrow & X\\
& \tilde{q}=(\tilde q_1,\tilde q_2,\alpha_1,\alpha_2) & \mapsto 
&\alpha_1\Phi_1(\tilde q_1)+\alpha_2\Phi_2(\tilde q_2). 
\end{array}\]
Since $T_{p_i}^m(X_i)$ is spanned by the derivatives of order less or equal 
than $m$ of $\Phi_i$ evaluated in $\Phi_i^{-1}(p_i)$ and 
$T_{p}^m(J(X_1,X_2))$ is spanned by the derivatives of order less or equal 
than $m$ of $\Phi$ evaluated in $\Phi^{-1}(p)$, the result follows.

\qed

It is clear that Theorem~\ref{maintheorem} is just the special case of 
Proposition~\ref{join} corresponding to $X_0 = \ldots = X_h = X$. 
Now we fix our attention on differential equations of second order
(the so-called Laplace equations) . In \cite{Terracini:12}
Terracini gives a characterization of varieties satisfying 
$$
T={K \choose 2}+m
$$
of such equations, where $K$ is the dimension of the variety and $m>0$
(a modern reference is \cite{Ilardi:99}, Theorem~1.5).
When $K$ grows up, $T$ becomes a huge number and we can ask, 
for instance, if secant varieties satisfy such a number of Laplace equations. 
We will see that the answer is negative, provided that we consider 
varieties without unexpected properties. Namely, let $X \subset \PP^r$ be 
a variety of dimension $n$ such that  $X$ is not $k$-defective for all $k$.
Thus, by Theorem~\ref{maintheorem} we have
\begin{equation}\label{dimen}
\dim(T_u^2(S^h(X)))\leq (h+1){{n+2}\choose n}-1.
\end{equation}
The previous inequality is strict when the $2$-osculating spaces at the 
generic points of $X$ intersect along subspaces of dimension greater than 
or equal to zero or they have dimension less than the expected one.
We assume that equality holds in (\ref{dimen}) and let  
$K:=\dim(\mbox{Sec}_h(X))=(h+1)n+h$. Thus, by Theorem \ref{maintheorem},  
$S^h(X)$ satisfies 
$$
T:={{K+2} \choose K}-(h+1){{n+2} \choose n} =
\frac{1}{2}n^2h^2+\frac{n^2h}{2}+nh^2+nh+\frac{1}{2}h^2+\frac{1}{2}h
$$  
Laplace equations. It is easy to see that we can write $T$ as
$$
T={K \choose 2}-\left(\frac{1}{2}n^2-\frac{1}{2}n-1\right)h-\frac{1}{2}n^2+
\frac{1}{2}n.
$$
Hence if $n\geq 2$ then  $T\leq {K \choose 2}$ and we are not in the 
hypothesis of \cite{Terracini:12}. If instead $n=1$ we have
$$
T \geq {K \choose 2} + h.
$$ 
Hence in this case we are in the hypothesis of \cite{Terracini:12}
and in the notation of \cite{Ilardi:99} we have:

\begin{Proposition}
Let $C$ be an integral non-degenerate curve in $\PP^r$ and suppose 
that $2h+1 \le r$, $h>0$. Then $S^h(C)$ is contained in a variety $U_q = 
\infty^t\PP^p$, such that the $\PP^q$ tangent at the points of a 
$\PP^p$ lie in a $\PP^{2K-t-m}$, with $0 \leq t \leq K-m$.
\end{Proposition}

Recall that an integral non-degenerate curve $C$ is not $k$-defective 
for any $k$ (in particular, $S^1(C)$ is a $3$-fold) and the osculating 
spaces to $C$ at a general point have always the expected dimension.
If $C$ is embedded as a rational normal curve 
in a projective space $\PP^r$ with $r$ big enough, 
then two osculating spaces $T^2_P(C)$ and $T^2_Q(C)$ do not intersect 
and $S^1(C)$ satisfies exactly four Laplace equations. We point out that
the case in which $X$ is a $3-$fold satisfying four Laplace equations is 
described in \cite{Terracini:12} by Terracini. He asserts that $X$ lies in 
$\PP^5$ or it is a $\infty^2$ of developable lines (with the tangent 
$\PP^3$ fixed along a line) or it lies in an ordinary developable $4-$fold.

\section{The Horace method and postulation \\of points in $\PP^3$}
\label{horace}

For any closed subscheme $Z \subset \PP^n$ and every integer $t$, let 
$$
\rho _{Z,t,n} : H^0({\bf {P}}^n,\mathcal {O}_
{{\bf {P}}^n}(t)) \to H^0(Z,\mathcal {O}_Z(t))
$$ 
be the restriction map. $Z$ is said to have maximal rank if
for every integer $t > 0$ the restriction map $\rho _{Z,t,n}$ has maximal rank, i.e. it is either injective or surjective.
Set $\rho _{Z,t}:= \rho _{Z,t,3}$.
Let $W$ be an algebraic scheme and $H \subset W$ an effective Cartier divisor of $W$. For any closed subscheme $A$ of $W$
let $\mbox{Res}_H(A)$ denote the residual scheme of $A$ with respect to $H$. 
Assume now that $W$ projective and fix $L\in \mbox{Pic}(W)$. By the very definition of residual
scheme we have the following exact sequence on $W$:
\begin{equation}\label{eq0}
0 \to \mathcal {I}_{\mbox{Res}_H(A)}\otimes L(-H) \to \mathcal {I}_A\otimes L \to \mathcal {I}_{A\cap H,H}\otimes (L_{\vert H}) \to 0
\end{equation}

From the long cohomology exact sequence of the exact sequence (\ref{eq0}) we get the following elementary form of the so-called Horace Lemma:

\begin{Lemma}\label{a0}
Let $W$ be a projective scheme, $H \subset W$ an effective Cartier divisor of $W$ and $A \subseteq W$ a closed subscheme of $W$.
Then:
\begin{itemize}
\item[(a)] $h^1(W,\mathcal {I}_A\otimes L) \le h^1(W,\mathcal {I}_{\mbox{Res}_H(A)}\otimes L(-H)) + h^1(H,
\mathcal {I}_{A\cap H,H}\otimes (L_{\vert H}))$;
\item[(b)] $h^1(W,\mathcal {I}_A\otimes L) \le h^1(W,\mathcal {I}_{\mbox{Res}_H(A)}\otimes L(-H)) + h^1(H,
\mathcal {I}_{A\cap H,H}\otimes (L_{\vert H}))$.
\end{itemize}
Now assume $h^1(W,L) = h^1(W,L(-H)) = h^1(H,L_{\vert H}) = 0$ and that $A$ is zero-dimensional.
If the restriction maps $H^0(W,L(-H)) \to H^0(A,L(-H)_{\vert \mbox{Res}_H(A)})$ and $H^0(H,L_{\vert H}) \to H^0(A,L_{\vert A\cap H})$
are surjective (resp. injective), then the restriction map $H^0(W,L) \to H^0(A,L_{\vert A})$ is surjective (resp. injective).
\end{Lemma}

\begin{Remark}\label{b1}
Fix integers $k \ge 3$, $a \ge 0$ and $b \ge 0$. Assume $(k,a,b)\notin 
\{(2,0,2),(3,2,0),(3,1,1),(4,0,5),(4,2,1),(4,2,0),(5,2,3),(6,5,0),(6,4,1)\}$. 
Let $Z \subset \PP^2$ be
a general union of $a$ triple points and $b$ double points. By the classification of all special linear systems of fat points
with multiplicity at most $3$ in the projective plane, we obtain that the restriction map $\rho _{Z,k,2}: H^0(\PP^2,\mathcal {O}_{{\bf {P}}^2}(k))
\to H^0(Z,\mathcal {O}_Z(k))$ has maximal rank.
\end{Remark}

For all integers $t \ge 2$ and $n \ge 2$ define the integers $a_{t,n}$ and $b_{t,n}$ by the relations:
\begin{equation}\label{eqa1}
{{n+2}\choose{2}}a_{t,n} + b_{t,n} = {{t+n}\choose{n}}, \ 0 \le b_{t,n} < 
{{n+2}\choose{2}}
\end{equation}
Set $a_t:= a_{t,3}$ and $b_t:= b_{t,3}$.

\begin{Notation}\label{a1}
For all integers $t \ge 2$, $a \ge 0$ and $b \ge 0$ such that
\begin{equation}\label{eqa2}
10a + 4b \le {{t+3}\choose{3}}
\end{equation}
we will say that the assertion $A_{t,a,b}$ is true if for a general union $Z \subset \PP^3$ of $a$ triple
points and $b$ double points the restriction map $\rho _{Z,t}$ is surjective.
\end{Notation} 

\begin{Notation}\label{a2}
For all integers $t \ge 2$, $a \ge 0$ and $b \ge 0$ such that
\begin{equation}\label{eqa3}
10a + 4b \ge {{t+3}\choose{3}}
\end{equation}
we will say that the assertion $B_{t,a,b}$ is true if for a general union $Z \subset \PP^3$ of $a$ triple
points and $b$ double points the restriction map $\rho _{Z,t}$ is injective.
\end{Notation}

\begin{Remark}\label{a3}
Let $Z \subset W$ be zero-dimensional schemes. Since $W$ is affine, the restriction map $H^0(W,\mathcal {O}_W) \to
H^0(Z,\mathcal {O}_Z)$ is surjective. Hence if $A_{t,a,b}$ is true, then $A_{t,x,y}$ is true for all pairs of
non-negative integers $(x,y)$ such that $x \le a$ and $x+y \le a+b$. Obviously, if $W \subset \PP^3$ and
$h^0(\PP^3,\mathcal {I}_Z(t))=0$, then $h^0(\PP^3,\mathcal {I}_W(t))=0$. Hence if $B_{t,a,b}$ is true,
then $B_{t,u,v}$ is true for all pairs $(u,v)$ of non-negative integers such that $u \ge a$ and $u+v \ge a+b$.
\end{Remark}

\begin{Remark}\label{a4}
Let $T$ be any integral projective variety and $V$ any linear system on $T$. For any integer $s > 0$ and general 
$P_1, \ldots, P_s \in T$ we have $\dim(V(-P_1 \ldots -P_s))= \min \{ 0, \dim(V) - s \}$.
\end{Remark}

\begin{Proposition}\label{a12.1}
Fix integers $k \ge 4$, $a \ge 0$ and $b\ge 0$. Set $\gamma := \max \{0,k-b-4\}$. If $10a+4b + 3\gamma +2k \le
{{k+3}\choose{3}}$, then $A_{k,a,b}$ is true.
\end{Proposition}

\proof
It is sufficient to show the existence of a disjoint union $E \subset \PP^3$ of $a$ triple points and $b$ double points
and ${{k+3}\choose{3}} - 10a-4b$ points such that $h^1(\PP^3,\mathcal {I}_A(k)) = 0$, i.e. such that
$h^0(\PP^3,\mathcal {I}_A(k)) = 0$. First assume $a \ge a_{k,2}$. Fix a plane $H \subset \PP^3$. Let $S$ be a general subset of $H$ with
$\mbox{card}(S) = a_{k,2}$ and let $E \subset \PP^3$ be the union of the triple points of $\PP^3$ such that
$E_{red} = S$. Hence $E\cap H$ is the union of $a_{k,2}$ triple points of $H$. By its very definition we have $0 \le b_{k,2}
\le 5$. First assume $0 \le b_{k,2} \le 2$. Let $M$ be a general subset of $H$ with cardinality $b_{k,2}$. By Remark \ref{b1}
we have $h^1(H,\mathcal {I}_{(E\cap H)\cup M}(k)) = h^0(H,\mathcal {I}_{(E\cap H)\cup M}(k)) = 0$. The residual scheme
$\mbox{Res}_H(E\cup M)$ is the union $E'$ of $a_{k,2}$ double points supported on $E_{red}$ and $b_{k,2}$ simple points.
Hence by Horace Lemma \ref{a0} and Remark \ref{a4} 
it is sufficient to show the existence of a union $B \subset \PP^3$ of $a-a_{k,2}$ triple points, 
$b$ double points and ${{k+3}\choose{3}}-10a-4b-b_{k,2}$ simple points such that $h^1(\PP^3,\mathcal {I}_B(k-1)) = 0$, i.e. such that
$h^0(\PP^3,\mathcal {I}_B(k-1)) = 0$. The idea is to use the same strategy with the integer $k-1$ instead of the integer
$k$. In order to do so, we need to look first at the case $3 \le b_{k,2} \le 5$. If $b=0$ (i.e. $\gamma = k-4$), then we do exactly the
same construction. Assume now $b>0$. Fix a general $P\in H$. Let $M' \subset H$ a general subset of $H$ with cardinality
$b_{k,2}-3$. Set $G:= E'\cup 2P\cup M'$. Hence $\mbox{length}(G\cap H) = 
{{k+2}\choose{2}}$. By Remark \ref{b1} we
have $h^0(H,\mathcal {I}_G(k)) = h^1(H,\mathcal {I}_G(k)) = 0$. We have $\mbox{Res}_H(G) = E'\cup \{P\}$. Hence we may try
to procede as in the case $b_{k,2} \le 2$ using $E'\cup \{P\}$ instead of $E'$ for the step concerning $\mathcal {O}_{\PP^3}
(k-1)$. We continue in one of the two ways according to the value of $b_{k-1,2}$ (of course, if $b \ge 2$) decreasing
each time the degree of the line bundle, until perhaps (say when we need to consider $\mathcal {O}_{\PP^3}(c)$)
we have only $a'< a_{c,2}$ triple points left. In this case we insert on $H$ the support of all $a'$ triple points
and use the integer $b_{c,2}+6(a_{k,2}-a')$ instead of the integer $b_{k,2}$. If $c = k$ we handle also
the case $a < a_{k,2}$ excluded at the beginning of the proof. The choice of the integer $\gamma$
and our assumption $10a+4b + 3\gamma +2k \le {{k+3}\choose{3}}$ shows that we can continue for
all lower degree line bundles until we exaust in in this way both the triple points and the double points.

\qed

The same proof (just changing the surjectivity part with the injectivity part of Horace Lemma \ref{a0}) gives the following result:

\begin{Proposition}\label{a12.2}
Fix integers $k \ge 4$, $a \ge 0$ and $b\ge 0$. Set $\gamma := \max \{0,k-b-4\}$. If $10a+4b \ge
{{k+3}\choose{3}}+ 3\gamma +2k$, then $B_{k,a,b}$ is true.
\end{Proposition}

\newpage

\vspace{0.5cm}
\noindent Edoardo Ballico \newline
Universit\`a degli Studi di Trento \newline
Dipartimento di Matematica \newline
Via Sommarive 14 \newline
38050 POVO (Trento) \newline
Italy \newline
e-mail: ballico@science.unitn.it

\vspace{0.5cm}
\noindent Cristiano Bocci \newline
Universit\`a degli Studi di Milano \newline
Dipartimento di Matematica ``F. Enriques'' \newline
Via Cesare Saldini 50 \newline
20133 MILANO \newline
Italy \newline
e-mail: cristiano.bocci@unimi.it
 
\vspace{0.5cm}
\noindent Enrico Carlini \newline
Universit\`a degli Studi di Pavia \newline
Dipartimento di Matematica ``F.Casorati'' \newline
Via Ferrata 1 \newline 
27100 PAVIA \newline
Italy \newline
e-mail: carlini@dimat.unipv.it

\vspace{0.5cm}
\noindent
Claudio Fontanari \newline
Universit\`a degli Studi di Trento \newline
Dipartimento di Matematica \newline
Via Sommarive 14 \newline
38050 POVO (Trento) \newline
Italy \newline
e-mail: fontanar@science.unitn.it

\end{document}